\documentclass[aos,dvips]{arximspdf}
\usepackage{graphics}

\doi{10.1214/09-AOS743}
\volume{38}
\issue{3}
\pubyear{2010}
\firstpage{1638}
\lastpage{1664}

\makeatletter
\newcommand{\X}{\mathbf{X}}
\newcommand{\x}{\mathbf{x}}

\newtheorem{lemma}{Lemma}[section]

\newtheorem{theorem}{Theorem}
\makeatother

\begin{document}
\begin{frontmatter}

\title{Successive normalization of rectangular arrays}
\runtitle{Successive normalization of rectangular arrays}

\begin{aug}
\author[A]{\fnms{Richard A.}
\snm{Olshen}\thanksref{tz}\ead[label=e1]{olshen@stanford.edu}\corref{}}\and
\author[B]{\fnms{Bala}
\snm{Rajaratnam}\thanksref{tm}\ead[label=e2]{brajarat@stanford.edu}}
\runauthor{R. A. Olshen and B. Rajaratnam}
\affiliation{Stanford University}
\address[A]{Departments of Health Research and Policy,\\
\quad Electrical Engineering, and Statistics \\
Stanford University\\
Stanford, California 94305-5405\\
USA\\
\printead{e1}} 
\address[B]{Department of Statistics\\
Stanford University\\
Stanford, California 94305-4065\\
USA\\
\printead{e2}}
\end{aug}
\thankstext{tz}{Supported in part by NIH MERIT Award R37EB02784.}
\thankstext{tm}{Supported in part by NSF Grants DMS-05-05303, DMS-09-06392, NIH Grant 5R01 EB001988-11REV and SUFSC2008-SUSHSTF2009-SMSCVISG200906.}

\received{\smonth{1} \syear{2009}}
\revised{\smonth{7} \syear{2009}}

%
\begin{abstract}
Standard statistical techniques often require transforming data to have
mean $0$ and
standard deviation $1$. Typically, this process of ``standardization''
or ``normalization''
is applied across subjects when each subject produces a single number.
High throughput genomic
and financial data often come as rectangular arrays where each
coordinate in one direction
concerns subjects who might have different status (case or control,
say), and each
coordinate in the other designates ``outcome'' for a specific feature,
for example, ``gene,''
``polymorphic site'' or some aspect of financial profile. It may
happen, when
analyzing data that arrive as a rectangular array, that one requires
BOTH the
subjects and the features to be ``on the same footing.'' Thus there may
be a
need to standardize across rows and columns of the rectangular matrix. There
arises the question as to how to achieve this double normalization. We propose
and investigate the convergence of what seems to us a natural approach
to successive
normalization which we learned from our colleague Bradley Efron. We
also study the
implementation of the method on simulated data and also on data that arose
from scientific experimentation.
\end{abstract}

%

\end{frontmatter}

\setcounter{footnote}{2}
\section{Introduction}

This paper is about a method for normalization, or regularization, of
large rectangular sets of numbers. In recent years many statistical
efforts have been directed towards inference on such rectangular
arrays. The exact geometry of the array matters little to the theory
that follows. Positive results apply to the situation where there are
at least three rows and at
least three columns. We explain difficulties that arise when either
numbers only two. Scenarios to which methodology studied here applies
tend to have many more rows than columns. Data can be from gene expression
microarrays, SNP (single nucleotide polymorphism) arrays, protein
arrays,
alternatively from large scale problems in imaging. Often there is
one column per subject with rows consisting of real numbers (as in
expression).
Subjects from whom data are
gathered may be ``afflicted,'' or not, with a condition that, while
heritable, is far
from Mendelian. A goal is to find rows, better groups of rows, by
which to distinguish afflicted subjects from other subjects. One can be
led to testing
many statistical hypotheses simultaneously, thereby separating rows
into those that are ``interesting'' for further follow-up and those
that seem not
to be. Genetic data tend to be analyzed by test ``genes'' (rows), beginning
with their being ``embedded'' in a chip, perhaps a bead. There may
follow a subsequent molecule that binds to the embedded
``gene''/molecule. A compound that makes use of the binding preferences
of nucleotides and to which some sort of ``dye'' is attached is then
``poured.'' The strength of binding depends upon affinity of the
``gene'' or attached molecule and the compound. Laser light is shined
on the object into which the test ``gene'' has been embedded; and from its
bending, the amount of bound compound is assessed from which the
amount of the ``gene'' is inferred. The basic idea is that a different afflicted
status may lead to different amounts of ``gene.''

With the cited formulation and ingenious technology, data may still
suffer from problems that have nothing to do with differences between
groups of subjects or with differences between ``genes'' or groups of them.
There may be differences in background---by column, or even by row.
Perhaps also ``primers'' (compounds) vary across columns for a given
row. For
whatever reasons, scales by row or column may vary in ways that do not
enable biological understanding. Variability across subjects could be
unrelated to afflicted status.

Think now of the common problem of comparing variables that can vary
in their affine scales. Because covariances are not scale-free, it
makes sense to compare in dimensionless coordinates that are centered
at 0, that
is, where values of each variable have respective means subtracted
off and are scaled by respective standard deviations. That way, each variable
is somehow ``on the same footing.''

Standardization, or normalization, studied here is done precisely so
that both ``subjects'' and ``genes'' are ``on the same footing.'' We
recognize one might require only that ``genes'' (or some ``genes'') be
on the same footing, and the same for ``subjects.'' The successive
transformations studied here apply when one lacks a priori
opinions that might limit goals. Thus, ``genes'' that result from the
standardization we study are transformed to have mean 0 and standard
deviation 1 across all
subjects while the same is true for subjects across all ``genes.'' How
to normalize? One approach is to begin with, say, row, though one could
as easily begin with columns. Subtract respective row means and divide
by respective
standard deviations. Now do the same operation on columns, then on
rows, and so on. Remarkably, this process tends to converge, even
rapidly in terms of numbers of iterations, and to a set of numbers that
have the described good limiting properties in terms of means and
standard deviations, by row
and by column.

In this paper we show by examples how the process works and demonstrate
for them that indeed it converges. We also include rigorous
mathematical arguments as to why convergence tends to occur. Readers
will see that the process and perhaps especially the mathematics that
underlies it are not as simple as we had hoped they would be. This
paper is only about convergence which is demonstrated to be
exponentially fast (or faster) for examples. The mathematics here does
not apply directly to ``rates.'' The Hausdorff dimension of the limit
set seems easy enough to study. Summaries will be reported elsewhere.

\section{Motivating example}

We introduce a motivating example to ground the problem that we address
in this paper. Consider a simple 3-by-3 matrix with entries generated
from a uniform distribution on [0, 1]. We standardize the initial matrix
$X^{(0)}$ by row and column, first subtracting the row mean from each
entry and then dividing each entry in a given row by its row standard
deviation. The matrix is then column standardized by subtracting the
column mean from each entry and then by dividing each entry by the
respective column standard deviation. In this section, these four steps
of row mean polishing, row standard deviation polishing, column mean
polishing and column standard deviation polishing require one iteration
in the process of attempting to row and column standardize the matrix.
After one such iteration, the same process is applied to resulting
matrix $X^{(1)}$ and the process repeated with the hope that successive
renormalization will eventually yield a row and column standardized
matrix. Hence these fours steps are repeated until ``convergence''
which we define as the difference in the Frobenius norm between two
consecutive iterations being less than $10^{-8}$.

\begin{figure}

\includegraphics{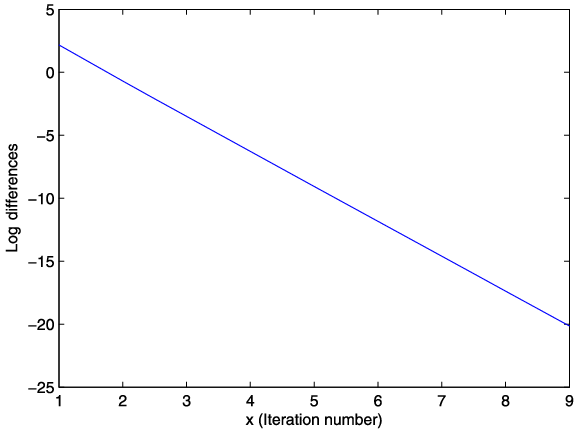}

\caption{Relative differences at each iteration on the log
scale---3-by-3 dimensional example.}\label{fig:21}
\end{figure}

In order to illustrate this numerically, we start with the following
3-by-3 matrix with independent entries generated from a uniform
distribution on [0, 1] and repeat the process described above:
%
\begin{equation}
X^{(0)} =
\left[\matrix{
0.1182 & 0.7069 & 0.4145 \cr
0.9884 & 0.9995 & 0.4648 \cr
0.5400 & 0.2878 & 0.7640
}\right]
.
\end{equation}
The successive normalization algorithm took 9 iterations to converge.
The initial matrix, the final solution and relative (and log relative)
difference for the 9 iterations are given below (see also Figure~\ref{fig:21}):
%
%
\begin{eqnarray}
X^{(\mathit{final})} &=&
\left[\matrix{
-1.2608 & 1.1852 & 0.0756 \cr
1.1852 & 0.0757 & -1.2608 \cr
0.0756 & -1.2608 & 1.1852
}\right]
,
%
%
%
\\
\hspace*{15pt}\mbox{Successive Difference} &=&
\left[\matrix{
\mbox{Iteration no.} & \mbox{Difference} & \log(\mbox{difference})\cr
1 & 8.7908 & 2.1737 \cr
2 & 0.5018 & -0.6895 \cr
3 & 0.0300 & -3.5057 \cr
4 & 0.0019 & -6.2862 \cr
5 & 0.0001 & -9.0607 \cr
6 & 0.0000 & -11.8337 \cr
7 & 0.0000 & -14.6064 \cr
8 & 0.0000 & -17.3790 \cr
9 & 0.0000 & -20.1516
}\right]
.
\end{eqnarray}

%

\begin{figure}

\includegraphics{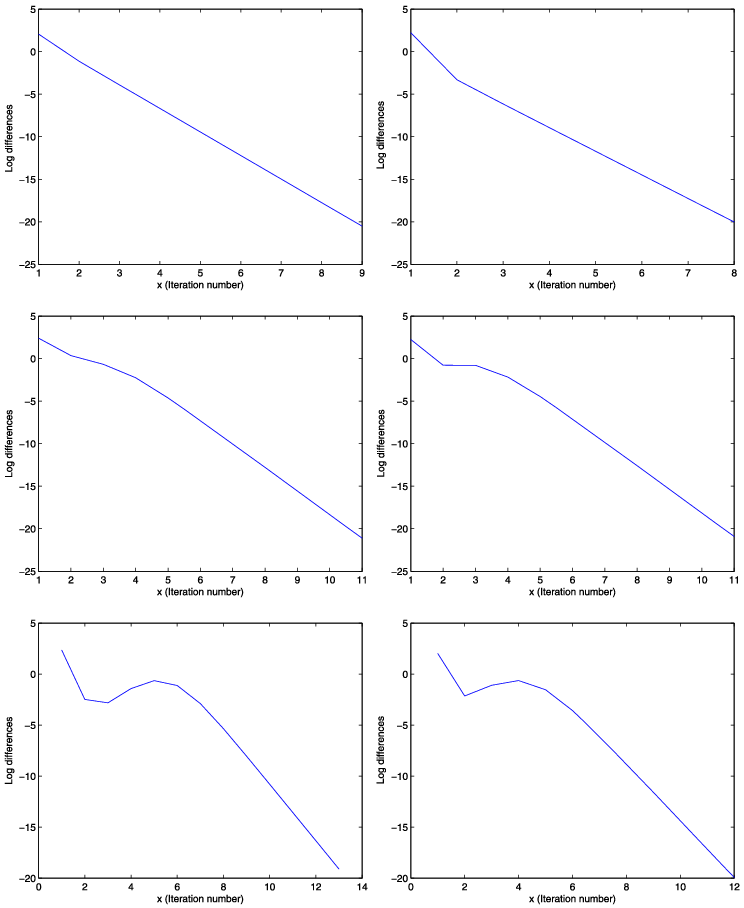}

\caption{Convergence patterns for the 3-by-3 example.}\label{fig:22}
\end{figure}

\noindent The whole procedure of 9 iterations takes less than 0.15
seconds on a standard modern laptop computer. We also note that the
final solution has effectively 3 distinct entries. When other random
starting values are used, we observe that convergence patterns can vary
in the sense that convergence may not be monotonic. The plots
in Figure~\ref{fig:22} capture the type of convergence patterns that
are observed in our simple 3-by-3 example.

%

 Despite the different convergence patterns that are observed,
when our successive renormalization is repeated with different starting
values, a surprising phenomenon surfaces. It always seems that the
process converges, and, moreover, the convergence is very rapid. One is
led naturally to ask whether this process will always converge and if
so under what conditions. These questions lay the foundation for the
work in this paper.

\section{Preliminaries}

We establish the notation that we will use by re-visiting a
normalization/standardization method that is traditional for
multivariate data. If the main goal of a normalization of a rectangular
array is achieving zero row and and column averages, then a natural
approach is to ``mean polish'' the row (i.e., subtract the row mean
from every entry of the rectangular array), followed by a column ``mean
polish.'' This cycle of successive row and column polishes is repeated
until the resulting rectangular array has zero row and and column
averages. The following theorem proves that this procedure attains a
double mean standardized rectangular array in one iteration where an
iteration is defined as constituting one row mean polish followed by
one column mean polish.
\begin{lemma}\label{lm:lemma1.3.1} Given an initial matrix ${\X}^{(0)}$, an
iterative procedure to cycle through repetitions of a row mean polish
followed by a column mean polish until convergence terminates in one step.
\end{lemma}

\begin{pf}
Let $\X^{(0)}$ be an $n\times k$ matrix, and define the following:%
\begin{eqnarray*}
\X^{(0)}&=& \bigl[X_{ij}^{(0)} \bigr],\\
\bar{X}_{i\cdot}^{(0)} &=&\frac{1}{k}\sum_{j=1}^{k} X_{ij}^{(0)}.
\end{eqnarray*}

Now the first part of the iteration, termed as a ``row mean polish,''
subtracts from each element its respective row mean:
\begin{eqnarray*}
\X^{(1)}= \bigl[X_{ij}^{(1)} \bigr]=X_{ij}^{(0)}-\bar{X}_{i\cdot}^{(0)}
.
\end{eqnarray*}

The second step of the iteration, termed a ``column mean polish,''
subtracts from each element of the current matrix its respective column
mean:
\[
\X^{(2)}= \bigl[X_{ij}^{(2)} \bigr]=X_{ij}^{(1)} - \bar{X}_{\cdot
j}^{(1)},
\]
where
\[
\bar{X}_{\cdot j}^{(1)}=\frac{1}{n}\sum_{i=1}^{n} X_{ij}^{(1)}.
\]

After the second step of the iteration it is clear that the columns sum
to zero; the previous operation enforces this. In order to prove that
the iterative procedure terminates at the second part of the iteration
it is sufficient to show that the rows of the current iterate sum to
zero. Now note that%
\begin{eqnarray*}
\X^{(2)}&=& \bigl[X_{ij}^{(2)} \bigr]\\
&=& \bigl[X_{ij}^{(1)} \bigr]-\bar{X}_{\cdot j}^{(1)}\\
&=&  \bigl(X_{ij}^{(0)}-\bar{X}_{i\cdot}^{(0)} \bigr) -  \Biggl(\frac{1}{n}\sum
_{r=1}^{n} X_{rj}^{(1)} \Biggr)\\
&=&  \bigl(X_{ij}^{(0)}-\bar{X}_{i\cdot}^{(0)} \bigr) -  \Biggl(\frac{1}{n}\sum
_{r=1}^{n} \bigl( X_{rj}^{(0)}-\bar{X}_{r\cdot}^{(0)} \bigr) \Biggr)
.
\end{eqnarray*}

It remains to show that the row sum of this matrix $\X^{(2)}$ expressed
as the elements of $\X^{(0)}$ sum to zero. So%
\begin{eqnarray*}
\sum_{j=1}^{k} X_{ij}^{(2)}
&=&\sum_{j=1}^{k}  \bigl(X_{ij}^{(0)}-\bar{X}_{i\cdot}^{(0)} \bigr)
- \sum_{j=1}^{k} \Biggl(\frac{1}{n}\sum_{r=1}^{n} \bigl( X_{rj}^{(0)}-\bar
{X}_{r\cdot}^{(0)} \bigr) \Biggr)\\
&=& \bigl(k\bar{X}_{i\cdot}^{(0)}- k\bar{X}_{i\cdot}^{(0)} \bigr)
- \frac{1}{n}\sum_{r=1}^{n}\sum_{j=1}^{k} \bigl( X_{rj}^{(0)}-\bar
{X}_{r\cdot}^{(0)} \bigr)\\
&=& \bigl(k\bar{X}_{i\cdot}^{(0)}- k\bar{X}_{i\cdot}^{(0)} \bigr) - \frac
{1}{n}\sum_{r=1}^{n} \bigl( k\bar{X}_{r\cdot}^{(0)} -k\bar
{X}_{r\cdot}^{(0)} \bigr)\\
&=&0-0\\
&=&0
.
\end{eqnarray*}

\noindent Note that the above double standardization is implicit in a
2-way ANOVA, and, though not explicitly stated, it can be deduced from
the work of Scheff\'{e} \cite{scheffe}. It is nevertheless presented
here, first, in order to introduce notation; second, as it is not
available in this form above in the ANOVA framework; and third, for the
intuition it gives since it is a natural precursor to the subject of
work in the remainder of this paper.
\end{pf}
%
\subsection{Example 1 (cont.): A 3-by-3 example with only mean polishing}
We proceed to illustrate the previous theorem on the motivating example
given after the introduction and draw contrasts between the two
approaches. As expected the successive normalization algorithm
terminates in one iteration. The initial matrix, the final solution and
the column and row standard deviations of the final matrix are given below:
%
%
\begin{eqnarray}
Y^{(0)} &=&
 \left[\matrix{
0.1182 & 0.7069 & 0.4145 \cr
0.9884 & 0.9995 & 0.4648 \cr
0.5400 & 0.2878 & 0.7640
} \right]
;
\\
%
%
Y^{(\mathit{column}\mbox{-}\mathit{polished})} &=&
 \left[\matrix{
-0.4307 & 0.0422 & -0.1333 \cr
0.4396 & 0.3347 & -0.0829 \cr
-0.0089 & -0.3769 & 0.2162
} \right]
;
\\
%
%
Y^{(\mathit{row}\mbox{-}\mathit{polished})}&=&Y^{(\mathit{final})}=
 \left[\matrix{
-0.2568 & 0.2161 & 0.0407 \cr
0.2091 & 0.1043 & -0.3134 \cr
0.0477 & -0.3204 & 0.2727
} \right]
;
\\
%
%
\mathit{Std}(\mathit{columns})&=&
 \left[\matrix{
0.1932 & 0.2311 & 0.2410
} \right]
;
\\
%
%
\mathit{Std}(\mathit{rows})&=&
 \left[\matrix{
0.1952 \cr
0.2257 \cr
0.2445
} \right]
.
\end{eqnarray}
We note that, unlike the motivating example, and as expected, the row
and column means are both $0$, but the standard deviations of the rows
and the columns are not identical, let alone identically $1$. Since
mean polishing has already been attained, and we additionally require
that row and column standard deviations to be $1$, it is rather
tempting to row and column standard deviation polish the terminal
matrix $Y^{(\mathit{final})}$ above. We conclude this example by observing the
simple fact that doing so results in the loss of the zero row and
column averages.
%
\subsection{The 2-by-2 problem} We now examine the successive row and
column mean and standard deviation polishing for a $2 \times2$ matrix
and hence illustrate that for the results in this paper to hold true,
the minimum of row $(k)$ and column dimension $(n)$ of the matrix under
consideration must be at least $3$, that is, $\min(k,n) \ge3$. Consider
the following general $2 \times2$ matrix:
\[
{\X^{(0)}} =
\pmatrix{ a &b \cr
c &d
}
.
\]

If $a < b$ and $c < d$,
then after one row normalization,
\[
{\X}^{(1)} =
\pmatrix{ -1 &1 \cr
-1 &1
}
;
\]
so the variance of its values, denoted by
$(\mathbf{S}^{(1)}_j)^2$, is 0.
Therefore, allowing for both inequalities to be reversed, and, assuming
that, for example, $a$, $b$, $c$ and $d$ are i.i.d. with continuous
distribution(s), then $P((\mathbf{S}^{(1)}_j)^2 = 0 ) = {1 / 2}$; in which case
the procedure is no
longer well defined.

A moment's reflection shows that if ${\X}$ is $n \times2$
with $n$ odd, then after each row normalization, each column has an odd
number of entries;
each entry being $-1$ or $+1$. However, each row has exactly one $-1$
and one $+1$.
Thus $(\mathbf{S}^{(m)}_i - \mathbf{S}^{(m+1)}_j)^2 \rightarrow0$ occurs only on a set
of product Lebesgue measure
0. With $\min(k, n) \ge3$, both tend to 1 a.e. as
$m \nearrow\infty$. Henceforth, assume $\min(k,n) \ge3$.

\section{Theoretical considerations}

For a matrix $\X$ as defined, take
$\x\in\mathcal{R}^{n \times k}$.
Let $\lambda$ denote Lebesgue measure on $\mathcal{R}^{n \times k}$ and
$P$ be the probability on the coordinates that renders them i.i.d. That is,
$\X_{ij}(\x_{ij}) = \x_{ij}$. $\X_{ij}$ are i.i.d.
$N(0, 1)$. What is more, for any orthogonal linear transformation
$\mathcal{O}$ that
preserves orthogonality, $\mathcal{O} \X\sim\X$, that is, $\mathcal{O} \X$
is distributed as $\X$ (see Anderson \cite{anderson} and Muirhead \cite{muirhead}).

Because $\lambda$ and $P$ are mutually absolutely continuous, if
$\mathcal{C} = \{$algorithm for successive row and column
normalization converges$\}$, then
$P(\mathcal{C}) = 1$ iff $\lambda(\mathcal{R}^{n \times k}
\setminus\mathcal{C})=0$, though only one direction is used.

For the remainder, assume that $P$ governs ${\X}$. Positive
results will be obtained for $3 \le \min(n,k) \le \max (n,k) < \infty
$. Our arguments require
notation, to which we turn our attention now. We redefine the notation
and symbols introduced in Lemma~\ref{lm:lemma1.3.1} as we now bring in row and column
standard deviation polishing. Define ${\X} = {\X^{(0)}}$:
\begin{eqnarray*}
{\X}^{(0)}_{i \cdot}
&=& {1\over k} \sum^k_{i=1} \X_{ij} ; \\
\bigl({\mathbf{S}}^{(0)}_i\bigr)^{2}
&=& {1\over k} \sum^k_{j=1} \bigl(\X_{ij} - {\bar{\X}}^{(0)}_{i \cdot
}\bigr)^2 \\
{}
&=& \Biggl({1\over k} \sum^k_{j=1} \X^2_{ij} \Biggr) -
\bigl({\bar{\X}}^{(0)}_{i \cdot}\bigr)^2;
\\
{\X}^{(1)} &=& \bigl[ \X_{ij}^{(1)} \bigr],\qquad \mbox{where }
\X^{(1)}_{ij} = {\bigl( \X_{ij} - {\bar{\X}}^{(0)}_{i \cdot}\bigr) / {\mathbf{S}
^{(0)}_i} } ;
\end{eqnarray*}
a.s.
$(\mathbf{S}^{(0)}_i)^2 > 0$ since $k \ge3 > 1$.

By analogy, set
\begin{eqnarray*}
\X^{(2)}
&=& \bigl[ \X_{ij}^{(2)} \bigr]\qquad \mbox{where } \X^{(2)}_{ij} =
{ \bigl(\X^{(1)}_{ij} - {\bar{\X}}_{\cdot j}\bigr) /\
{\mathbf{S}^{(1)}_j} } ; \\
{\bar{\X}}^{(1)}_{\cdot j}
&=& {1\over n} \sum^n_{i=1} \X^{(1)}_{ij} ; \\
\bigl(\mathbf{S}^{(1)}_j\bigr)^2
&=& {1\over n} \sum^n_{i=1} \bigl(\X^{(1)}_{ij} - {\bar{\X}}^{(1)}_{\cdot j}\bigr)^2 .
\end{eqnarray*}
Arguments sketched in what follows show that a.s.
$(\mathbf{S}^{(1)}_j)^2 > 0$ since $n \ge3$.

For $m$ odd, $\X^{(m)}_{ij} =
{( \X^{(m-1)}_{ij} - {\bar{\X}}^{(m-1)}_{i \cdot}) /\
{\mathbf{S}^{(m-1)}_i} } $
with ${\bar{\X}}^{(m-1)}_{i \cdot}$ and
$(\mathbf{S}^{(m-1)}_i)^2$
defined by analogy to previous definitions.

For $m$ even, $\X^{(m)}_{ij} =
{( \X^{(m-1)}_{ij} - {\bar{\X}}^{(m-1)}_{\cdot j}) /\
{\mathbf{S}^{(m-1)}_j} }$,
again with ${\bar{\X}}^{(m-1)}_{\cdot j}$ and
$\mathbf{S}^{(m-1)}_j$
defined by analogy to earlier definitions.

\subsection*{Back to the general problem}

We first note that because the process we study is a coordinate
process, there is no difference between regular conditional
probabilities and regular conditional distributions (see Durrett \cite{durrett}, Section~4.1.c, pages~33 and~229--331 for more details). They
can be computed as densities with respect to Lebesgue measure on a
finite Cartesian product${}\times \operatorname{Sph}(q)$ where $q=k$ refers to after
row normalization and where $q=n$ if subsequent to column
normalization. In a slight abuse of notation, for any positive integer
$q$ we define $\operatorname{Sph}(q) = \{ \x\in\mathcal{R}^q\dvtx  \| \x\|^2 = q\}$ where
the norm is the usual Euclidean norm.

Let $\{r_{ij} \dvtx  i = 1, \ldots, n; j=1, \ldots, k\}$ be a set of $n\cdot k$
rational numbers, not all 0. Obviously,
\[
P\biggl(\bigcup_{r_1, \ldots, r_{nk}} \sum_{i,j} r_{ij} \X_{ij} = 0\biggr) = 0 .
\]
An inductive argument involving conditional densities shows that
%
\begin{equation}\label{lincombzero}
P\Biggl(\bigcup^\infty_{m=1} \bigcup_{r_1, \ldots, r_{nk}} \sum_{i,j} r_{ij}
\X^{(m)}_{ij} = 0\Biggr) = 0.
\end{equation}
Consequently
\begin{eqnarray*}
P\Biggl(\Biggl(\bigcap^\infty_{m-1} \bigcap^n_{i=1} \bigl(\mathbf{S}^{(2m)}_i \bigr)^2 > 0\Biggr)
\cap\Biggl(\bigcap^\infty_{m=1}
\bigcap^k_{j=1} \bigl(\mathbf{S}^{(2m-1)}_j \bigr)^2 > 0\Biggr) \Biggr) = 1 .
\end{eqnarray*}

\noindent
Further, a.s. $\X^{(m)}$ is defined and finite for every $m$.

What we know that the $t$ distribution (Efron \cite{efron}) and
geometric arguments require that
${\X}^{(1)}$ can be viewed as having probability distribution
on
$\mathcal{R}^{n\times k}$
that is the $n$-fold product of independent uniform
distributions on $\operatorname{Sph}(k) \times \cdots \times \operatorname{Sph}(k)$.
For the sake of clarity we note again that conditional probabilities
are always taken to be ``regular'' and
concentrated on the relevant product of spheres. Readers will note that
in the cited
arguments for
$\sum_{i,j} r_{ij} \X^{(m)}_{ij}$, $(\mathbf{S}^{(2m)}_i)^2$ and
$(\mathbf{S}^{(2m-1)}_j)^2$, relevant
conditional probabilities and densities are used. Of course, after
the $m$th row standardization
$\X^{(2m-1)} = [\X^{(2m-1)}_{ij} \dvtx  i=1, \ldots, n; j=1, \ldots, k]$, and
analogously for column standardization and $\X^{(2m)}$.

As an aside, write $g_1(\X) = \X^{(1)}$ on
$\{\min (\mathbf{S}^{(0)}_i )^2 > 0\}$. Elsewhere, let
$g_1 (\X)=0$. Necessarily $g_1$ depends on $\X$ only through its
direction in
$\mathcal{R}^{n\times k}$. Equivalently, $g_1$ is a function of $\X$ that is
homogeneous of degree $0$. Moreover,
$g_1 (\X) \sim g_1 (\mathcal{O} \X)$ for all orthogonal linear $\mathcal{O}$
on $\mathcal{R}^{n\times k}$. $\X^{(1)}$ has independent rows, each uniformly
distributed on $\operatorname{Sph} (k)$.

We turn now to study $\X^{(2m-1)}_{ij}$ as $m$ increases without
bound.
Note first that for $m=1, 2, \ldots, \X^{(2m-1)}$ has joint distribution that is unchanged if
two columns of $\X$ are transposed, therefore if two columns of $\X
^{(2m-1)}$ are
transposed. Since every permutation is a product of transpositions,
$\X^{(2m-1)}$ is column exchangeable. Each is also row exchangeable.

Write $\pi$ for a permutation of the integers $\{ 1, \ldots, k\}$; let
$\Pi$ be the finite
$\sigma$-field of all subsets of $\{ \pi\}$. The marginal
probability induced on $\{ \pi\}$ from the joint distribution of
$(\X, \{ \pi\})$ is discrete and uniform, assigning probability $1/k
!$ to each~$\pi$.

\noindent Write $\mathcal{G}^{(i)}_{2m-1}$ to be the $\sigma$-field,
\[
\mathcal{F}\bigl(\bigl[ \bigl(\X^{(q)}_{ij}\bigl)^2 \dvtx  j=1, \ldots, k; q = 2m -1, 2m+1, \ldots\bigr]\bigr)
\times\Pi.
\]
\renewcommand{\thetheorem}{4.1}
\begin{theorem}
$E\{ (\X^{(1)}_{i\pi(1)})^2 | \mathcal{G}_{2m-1} \} =
(\X^{(2m-1)}_{i\pi(1)} )^2$ a.s. for $m=1, 2, \ldots.$
\end{theorem}

\begin{pf}
Write
$(\X^{(2m-1)}_{i\pi(1)})^2 = \sum^k_{l=1} (\X^{(2m-1)}_{il})^2
I_{[\pi(1)=l]}$ where
$I_A$ is the indicator function of the event A. Obviously,
$(\X_{i\pi(1)}^{(2m-1)})^2$ is
$\mathcal{G}_{2m-1}^{(i)}$ measurable;
$\{\alpha_{q_1, q_2}\}$ is a set of real numbers, doubly indexed by
$\{\alpha_{q_1, q_2}) \in F$; $F$ is a finite subset of
$\{1, \ldots, k\} \times\{ 1, 2, \ldots \}$. Form
$B = \{( \X^{(2m-1+q_2)}_{iq_1})^2 \le\alpha_{q_1, q_2}; (q_1, q_2) \in
F\}$.

Note that $\mathcal{G}_{(2m-1)}^{(i)}$ is generated by $\{B \times Q\}$,
$B$ of the cited
form and $Q \in\Pi$. In particular, each $B \times Q \in\mathcal{
G}_{2m-1}^{(i)}$. Proof of our claim is complete if we show that for
$m=2, 3, \ldots,$
\begin{eqnarray*}
\int_{B\times Q} \bigl(\X^{(2m-1)}_{i\pi(1)}\bigr)^2 =
\int_{B\times Q} \bigl(\X^{(1)}_{i\pi(1)}\bigr)^2 .
\end{eqnarray*}

The left-hand side of the display can be expressed as
\begin{eqnarray*}
E\Biggl\{ \sum^k_{l=1} \bigl(\X^{(2m-1)}_{il}\bigr)^2 I_{[\pi(1)=l]} I_{[\pi
\in Q]} I_B \Biggr\}
=E\Biggl\{ I_B \sum^k_{l=1} \bigl(\X^{(2m-1)}_{il}\bigr)^2 I_{[\pi(1)=l]} I_{\pi\in Q]}\Biggr\}.
\end{eqnarray*}

Now, for any $\pi$, the expression inside the sum is $k$ if
$\pi\in Q$ and~$0$ if not.
%
%
%
So, the expectation factors into $P(B) P(Q).$
Retracing steps shows clearly that
$(\X^{(2m-1)}_{il})^2$ in the computation just completed could be
replaced by
$(\X^{(1)}_{il})^2$ with  equalities remaining true. The claim is
now proven.
\end{pf}

\subsection*{Convergence of $(\X^{(m)}_{ij})^2$}

The backwards martingale convergence theorem (see Doob \cite{doob})
entails that
$(\X^{(2m-1)}_{i\pi(1)})^2$ converges a.s. as $m \rightarrow\infty$.
So, for each fixed $j\in\{1, \ldots, k\}$, $(\X^{(2m-1)}_{i\pi(1)})^2I_{[\pi(1)=j]}$
converges a.s. It follows that $[(\X^{(2m-1)}_{ij})^2]$ converges a.s.
as $m\rightarrow\infty$.

If previous arguments are perturbed so that $\pi$ denotes a
permutation of
$\{1, \ldots, n\}$ with $(\X^{(1)}_{i \pi(1)})^2$ replaced by\vspace*{1pt}
$(\X^{(2)}_{\pi(1) j})^2, \mathcal{G}^{(i)}_{2m-1}$ by
$\mathcal{G}^{(i)}_{2m}$,
$(\X^{(2m-1)}_{i\pi(1)})^2$ by $(\X^{(2m)}_{\pi(1) j})^2$ and
$\sum^k_{l=1} (\X^{(2m-1)}_{il})^2$ by
$\sum^n_{l=1} (\X^{(2m)}_{lj})^2$, then one concludes\vspace*{1pt} that also
$[(\X^{(2m)}_{ij})^2]$ converges a.s. as $m\rightarrow\infty$.
Without further argument it is unclear if the a.s. limits along odd,
respectively even, indices are the same; and it is crucial to what
remains that this is
in fact true.

Obviously,
$\bigcap^\infty_{m=1} \mathcal{G}^{(i)}_{2m-1} =
\bigcap^\infty_{m=1} \mathcal{G}^{(i)}_{2m}$, so in a certain sense
measurability is the same. Obviously, too, randomization of index is by
columns in the
first case and by rows in the second. But now a path to the required
conclusion presents
itself. Given our success in proving a.s. convergence along odd indices
after randomizing columns and along even indices randomizing rows, and given
that a requirement of our approach is that these two limits be identical
a.s., perhaps there is a path which would allow for the simultaneous
randomizing of both columns and
rows. Fortunately, that is the case. Thus, let
$\pi_1$ be a permutation of $\{1, \ldots, n\}$ and $\pi_2$ be a
permutation of $\{1, \ldots, k\}$. With obvious product formulation of
governing probability
mechanism and further obvious formulation of decreasing $\sigma
$-fields, as an example of
what can be proved,
\begin{eqnarray*}
E\bigl(\bigl( \X^{(1)}_{\pi_2 (1) \pi_2 (1)} \bigr)^2 | \mathcal{G}_2\bigr) =
\bigl(\X^{(2)}_{\pi_1 (1) \pi_2(1)}\bigr)^2 \qquad\mbox{a.s.}
\end{eqnarray*}
\noindent
From this arguments for the display there are several paths by
which one concludes that a.s., simultaneously for all $(i,j), [(\X
^{(m)}_{ij})^2]$ converges. Dominated convergence requires that the
limit random matrix has expectation 1 in all coordinates. As a
consequence of this convergence, a.s. and simultaneously for all
$(i,j), [(\X^{(2m+1)}_{ij})^2] - [(\X^{(2m)}_{ij})^2] \rightarrow0$
as $m\rightarrow
\infty$.

\subsection*{Convergence of $[\X_{ij}^{(m)}]$}

We turn now to key ideas in extending our argument\vspace*{-1pt} that $[(\X
_{ij}^{(m)})^2]$ converges almost surely simultaneously to the same
conclusion with the square removed. To limit notational complexity, we
study first only odd indices as m grows without bound. Conclusions are
identical for even indices, and by extension for indices not
constrained to be odd or even.

A first necessary step is to show that for arbitrary $j$,
\[
P\Bigl\{\overline{\lim_m} \bigl(S^{(2m+1)}_j \bigr)^2 > 0\Bigr\}=1.
\]

To that end, let A be the event $[ \overline{\lim}_m
(S^{(2m-1)}_j )^2 = 0]$. Obviously,\vspace*{-2pt} $A =\{\x\dvtx  (S^{(2m-1)}_j )^2
\rightarrow0\} = \{\x\dvtx  S^{(2m-1)}_j \rightarrow0\}$ regardless of
square roots taken. We show that $P\{A\}=0$.

By way of contradiction, suppose that $P\{A\}> 0$. Write
\[
\bigl(S^{(2m-1)}_j\bigr)^2 = {1\over n} \sum^n_{l=1} \bigl(\X^{(2m-1)}_{lj}\bigr)^2 - \bigl({\bar
{\X}}^{(2m-1)}_{\cdot j}\bigr)^2.
\]

We know that the first term tends to 1 a.s. on $A$. Therefore, also the
second term tends to 1 a.s. on $A$. Since for $m>1$, $ \X^{(2m-1)}_{lj}
$ is bounded \vspace*{-2pt}a.s., $\overline{\lim}_m   \max_l   (\X
^{(2m-1)}_{lj})(\x)$ is a finite-valued random variable $C=C(\x)$.
Simple considerations show that the only possibilities are that for all
$l$, $C=+1$ or $C=-1$ a.s. Similar arguments\vspace*{1.5pt} show that
$\mathop{\underline{\lim}}_m \min_l (\X^{(2m-1)}_{lj})(\x)$ is $+1$ or
$-1$. It is clear that there is no sequence $\{m_k \}$ of $\{m\}$ along which

\[
-1 < \mathop{\underline{\lim}}_m   \min_l   \bigl(\X^{(2m_{k}-1)}_{lj}\bigr)(\x)
\leq \overline{\lim_m}   \max_l   \bigl(\X^{(2m_{k}-1)}_{lj}\bigr)(\x)
< 1.
\]

It follows that $\lim_m   \X^{(2m-1)}_{lj}$ exists a.s. on $A$
and that the limit of the sequence is $+$1 or $-$1 on $A$.

Recall that $\X\sim-\X$, and this equality is inherited by all joint
distributions of $\X^{(m)}$. However, when $\X$ is replaced by its
negative, any limit of $\X^{(2m-1)}_{lj}$ becomes its negative; a.s.
convergence is a property of the probability distribution of $\X$, and
$\{\x\dvtx  \X^{(2m-1)}_{lj} \mbox{ converges} \}= \{\x\dvtx  -\X^{(2m-1)}_{lj}
\mbox{ converges} \}$. Therefore, the only possibility is that
$(\overline{\X}^{(2m-1)}_{\cdot j} )^2 \rightarrow0 $ and
$(S^{(2m-1)}_{j} )^2 \rightarrow1$. The upshot is that $P\{A\}=0$; $
P\{\overline{\lim}_m (S^{(2m-1)}_j )^2 > 0\}=1 $.

Again, let us fix $j$. Consider a sample path of
$\{\X^{(2m_q-1)}\}$ along which
$\lim_{m} (S_j^{(2m_q-1)})^2 = D >0$. Clearly,
$\{i\dvtx  {\overline{\lim}}_{m_q} | \X^{(2m_q-1)}_{ij} | > 0\} \not=
\varnothing$. Indeed,\vspace*{-1pt} let
$E=E(j) = \{i\dvtx \mbox{for some }
\{m_q\} = \{m_q (i)\}, {\overline{\lim}}_{{m_q}(i)} |
\X_{ij}^{(2m_q -1)} | > 0\}$.
Row and column exchangeability of $\X^{(m)}$ necessarily means that
the cardinality of $E$ is at least~2.

Let $i_0 \not= i_1 \in E$. Because $\min(n, k) \ge3$,
there is a further subsequence of $\{\{ m_q (i) \}\}$---again, for
simplicity, write it as $\{m_q\}$---along which
\begin{eqnarray*}
&&\lim_{m_q} \big| \X^{(2m_q - 1)}_{{i_0} j} \big| \quad\mbox{and}\quad
\lim_{m_q} \big| \X^{(2m_q)}_{{i_0} j} \big| \qquad\mbox{both exist};
\\
&&\lim_{m_q} \big| \X^{(2m_q - 1)}_{{i_1} j} \big|\quad \mbox{and}\quad
\lim_{m_q} \big| \X^{(2m_q )}_{{i_1} j} \big|\qquad\mbox{ both exist and }
\\
&&\lim_{m_q} \big| \X^{(2m_q)}_{{i_0} j} - \X^{(2m_q - 1)}_{i_{1}j} \big|
\qquad\mbox{exists and is positive.}
\end{eqnarray*}

The first two requirements can always be met off of the set of
probability 0 implicit in (\ref{lincombzero}). That the third can be
met as well is a consequence of the argument just concluded. In any
case, if there were no such subsequence, then our proof would be
complete because all $\X^{(m)} _{ij}$ for $j$ fixed tend to the same
number. But now, write
\begin{eqnarray*}
&&\bigl(\X^{(2m_q)}_{{i_0}j} - \X^{(2m_q - 1)}_{{i_0}j}\bigr) -
\bigl(\X^{(2m_q)}_{{i_1}j} - \X^{(2m_q - 1)}_{{i_1}j}\bigr)
\\
&&\qquad=\bigl(\X^{(2m_q - 1)}_{{i_0}j} - \X^{(2m_q - 1)}_{{i_1}j}\bigr)
\bigl(S^{(2m_q - 1)}_j - 1\bigr) / S^{(2m_q - 1)}_j .
\end{eqnarray*}
\noindent Since\vspace*{1.5pt}
$(\X^{(2m_q)}_{{i_0} j} )^2 - (\X^{(2m_q - 1)}_{{i_0} j} )^2
\rightarrow0$ a.s.,
and likewise with $i_0$ replaced by $i_1$, the first expression of the
immediately previous
display has limit $0$. Thus, so too does the second expression. This
is possible only if
$S^{(2m_q -1)}_j \rightarrow1$ (where we have taken the positive
square root). Further,
\begin{eqnarray*}
\X^{(2m_q)}_{{i_0} j} - \X^{(2m_q - 1)}_{{i_0} j} =
{{ \X^{(2m_q - 1)}_{{i_0} j} (S^{(2m_q - 1)}_j - 1) +
{\bar{\X}}^{(2m_q - 1)}_{\cdot j} } \over
{S^{(2m_q - 1)}_j} }
.
\end{eqnarray*}

\noindent
As a corollary to the above, one sees now that
${\bar{\X}}^{(2m_q - 1)}_{\cdot j} \rightarrow0$.
Since the original $\{m_q\}$ could be taken to be an arbitrary
subsequence of
$\{m\}$, we conclude that:
\begin{itemize}
\item$S^{(2m_q -1)}_j \rightarrow1$ a.s.,
\item${\bar{\X}}^{(2m_q -1)}_{\cdot j} \rightarrow0$ a.s. and
\item$\X^{(m)}_{ij}$ converges a.s.
\end{itemize}

\noindent Now replace arguments for (i) and (ii) on columns by
analogous arguments on rows. Deduce that every infinite subsequence of
positive integers has a subsequence along which our desired conclusion obtains.

\section{Properties of successive normalization}

We now comment on theoretical properties of successive normalization.
In particular, we elaborate on the generality of the result by showing
that the Gaussian assumption is not necessary and serves only as a
convenient choice of measure. We also discuss convergence in Lebesgue
measure and the domains of attraction of successive normalization.

\subsection{Choice of probability measure}

Write $\lambda$ for Lebesgue measure on $\mathcal{R}^{n\times k}$. Thus
$\mathbf{x} \in\mathcal{R}^{n\times k}$ is an
$n \times k$ rectangular array of real numbers. Let $P$ be a measure on
the Borel sets of
$\mathcal{B}$ of $\mathcal{R}^{n\times k}$ that is mutually absolutely
continuous with respect of $\lambda$. By this we mean that if
$B \subset\mathcal{R}^{n\times k}$ is Borel measurable, then
$\lambda(B) = 0$ iff $P(B) = 0$. One obvious example of such a $P$
is the measure on $\mathcal{B}$ that renders its $nk$ coordinates i.i.d.
Gaussian. If
$x_{rc}$ is the $(r, c)$ coordinate of
$\mathbf{X} \in\mathcal{R}^{n\times k}$, and $P^{(r,c)}$ is the marginal
measure of $P$ on its $(r, c)$ real coordinate, then
$P^{(r,c)} ((-\infty, x_0]) = \Phi(x_0)$ where $\Phi$ is the
standard normal cumulative distribution function. This is what we
mean by $P$ henceforth. The i.i.d. specification requires that $P$
is an $nk$-fold product of identical probabilities, and because $P$
is now defined uniquely for all product rectangles, it is
defined uniquely for all $B \in\mathcal{B}$. It is obvious that $\Phi$
being mutually absolutely continuous with respect to one-dimensional
Lebesgue measure means that the product measure $P$ is mutually
absolutely continuous with respect to the $nk$-fold Borel product
measure $\lambda$.

Now, let $f_1, f_2, \ldots$ be a sequence of $\mathcal{B}$-measurable functions,
$\mathcal{R}^{n\times k} \rightarrow\mathcal{R}^{n\times k}$. Then
$\{f_m$ converges$\} = \{ \overline{\lim} f_m = \underline{\lim}
f_m$ simultaneously for all $nk$ coordinates\}. When $f_m (\mathbf{x})$ converges, then
$\lim f_m (\mathbf{x}) = (\overline{\lim} f_m (x)) I_{\{f_m \mathrm{\ converges}\}}$, where for any set $C, I_C$ is its indicator. Because
each $f_m$ is $\mathcal{B}$-measurable,
$\lim f_m (\mathbf{x})$ is also $\mathcal{B}$-measurable. If we are given
that $\{f_m\}$ converges $P$-almost everywhere, then $\{f_m $
converges\} is a Borel set of
$P$-measure 1. Its complement has $P$-measure $0$ and,
therefore, $\lambda$ measure $0$. It follows that $\{f_m\}$ converges
except for a Borel subset of
$\mathcal{R}^{n\times k}$ of Lebesgue measure $0$.

In the present paper, the $(r, c)$ coordinate of $\{f_m\}$ is the set
of successive normalizations of the initial real entry multiplied by
the indicator of the subset of $\mathcal{R}^{n\times k}$ that is
\{successive normalizations possible\}. The latter is clearly a
Borel subset of $\mathcal{R}^{n\times k}$ of $P$-measure 1, so its
complement is a Borel set of $\lambda$-measure $0$. It is immaterial
for convergence
whether the first initialization is by row or by column. Readers
should note that any study of asymptotic properties of $\{f_m\}$
under $P$ may show that $\{f_m (\mathbf{x})\}$ has properties such as
row and column exchangeability, where the coordinate functions are
taken to be random variables. They should note, too, that: (i)
changing the original measure to one more conducive to computation is
standard, and is what happens with ``exponential tilting'' as it
applies to the study of large deviations of sums; (ii) the
interplay between measure and topology, as is utilized here, is a
standard approach in probability theory that is applied seldom to
statistical arguments; the lack thus far owes only to necessity.

\subsection{Convergence in $p$th mean for Lebesgue measure}
Whenever standardization is possible, after one
standardization the sum over all $nk$ coordinates
of squares of respective values is bounded by $nk$, it follows from
dominated convergence that
as $m$ grows without bound, each term converges
not only $P$-almost everywhere but also in $p$th mean for every
$\infty>p \ge1$ so long as $P$ remains the
applicable measure. Because $\lambda$ is not a
finite measure, convergence in $p$th mean for underlying measure
$\lambda$ does not follow. It is impossible for
such convergence to apply to Lebesgue measure on
the full Euclidean space $\mathcal{R}^{n\times
k}$. This is because there is a set of positive
Lebesgue measure whose members are not fixed
points for normalization by both rows and
columns. No matter which normalization comes
first in any infinite, alternating sequence, the
normalization is invariant to fixed scale
multiples of each $\mathbf{x} \in\mathcal{R}^{n\times
k}$, and, in particular, to fixed scale multiples of
$\mathbf{x}$ in the set of positive Lebesgue measure
just described. Obviously, no further arguments
are required; and convergence in $p$th mean is immediate if Lebesgue measure
$\lambda$ is restricted to a bounded subset of $\mathcal{R}^{n\times k}$.

\subsection{Domains of attraction} \label{domains_of_attraction}

Reviewers of research presented here have
wondered if we can describe simply what
successive and alternating normalization does to
rectangular arrays of data, beyond introductory
comments about putting rows and columns on an
equal footing and the analogy to computing
correlation from covariance. We begin our reply
here though details await further research and a
subsequent paper. Please remember invariance to
either row or column normalization (when possible) to scale multiples of
$\mathbf{x} \in\mathcal{R}^{n\times k}$. In other
words, results of normalization are constant
along rays defined by these multiples, and
without loss of generality we can assume that $\mathbf{x}$
lies in an $n$-fold product of $\operatorname{Sph}(k)$ where
each of the $n$ components is orthogonal to the
linear one-dimensional subspace of $\mathcal{R}^k$ consisting of its
equiangular line; call it $\operatorname{Sph}(k)\setminus
\{\mathbf{1}\}$. Thus, without loss we assume that
the object under study is $\mathbf{X}^{(m)}$ subject
to a row normalization of the process of
successive normalization. We study only subsets
of the $n$-fold product space that was described and that is
complementary to
$\{ \bigcup^\infty_{m=1} \bigcup_{r_1, \ldots, r_k}
\sum_{i,j} r_{i,j} \mathbf{X}^{(m)} = 0\}$
where $\{ r_{ij}\}$ are rational. The set in
braces just before the comma has $P$-measure $0$, and we know that on
its complement, $X^{(m)}$ can be defined for all $m=1,2,\ldots.$

Because normalization always involves subtraction
of a mean and division by a standard deviation
and because each $\mathbf{X}^{(m)}$ is row and
column exchangeable, the limiting process we
study here when $P$ applies seems, at first
glance, to be analogous to ``domains of
attraction'' as that notion applies to sequences
of i.i.d. random variables. One obvious
difference is that here limits are almost sure, unlike distribution.
While a.s. limits
of $\mathbf{X}^{(m)}$ are shown to have row and
column means $0$ and row and column standard
deviations $1, n\times k$ arrays of real numbers
with this property are obviously the only fixed
points of the alternating process studied
here. The Hausdorff dimension of the set of
fixed points is not difficult to compute, but we
have been unable thus far to give rigorously
supported conditions for the domain of attraction
(in the sense described) of each fixed
point. The simple case for which domains of
attraction for limits in distribution were
described was a major development in the history
of probability (see Feller \cite{feller}, Gnedenko and Kolmogorov \cite{gnedenkokolmogorov}, Zolotarev \cite{zolotarev}). We report some
intuitive results, and next we look at a mathematical
question that arose in our study of domains of
attraction for which at present we have only a heuristic argument.

Is there a set $E\subset\mathcal{R}^{n\times k}$ for which
$P(E)=1$; $X^{(m)} (\mathbf{x})$ converges for $\mathbf{x} \in E$, and for each fixed
$i   \lim\mathbf{X}_{ij} (\mathbf{x}) \neq\lim\mathbf{X}_{ij'} (\mathbf{x})$
for all
$j\neq j'$ if $\mathbf{x}\in E$? Clearly one could
ask the equivalent question for each fixed $j$, and a corresponding subset
$E'\subset\mathcal{R}^{n\times k}$. We conjecture
the existence of $E \cap E'$. However, arguments available thus far do
not confirm existence rigorously. Therefore,
suggestions regarding domains of attraction as
$\mathbf{X}^{(m)}$ grows without bound should be taken as only heuristic
for now.

Given that $( \mathbf{S}_i^{(m)})^2 \rightarrow1$ on a subset of
$\mathcal{R}^{n\times k}$ with complementary
$P$-measure~$0$, therefore Lebesgue measure $0$,
almost surely ultimately [meaning for $m=m(\mathbf{x})$ large enough],
row normalization does not
perturb the (strict) sort of any row. By
analogy, almost surely ultimately column normalization
does not perturb the (strict) sort of any
column. Thus, a.s. ultimately, successive and
alternative row and column normalization do not
perturb any row or column sort. This joint
strict sort determines an open subset of $\mathcal{R}^{n\times k}$. If we restrict ourselves to
rows, then we may speak more precisely of
$n$-fold products of sorts of members of
$\operatorname{Sph}(k)\setminus\{ \mathbf{1}\}$. For a fixed row
there are $k!$ strict sorts of the squared
entries. Given a fixed sort of the squares there are, say, $f(k)$
assignments of sign to the square roots of the $k$ entries so that the
sum of entries for that row is $0$. For any assignments of signs to the
necessarily nontrivial values of the $k$ squares that renders the sum
of entries for the row $0$, there is always a scaling so that the sum
of squares is any fixed value so that the variance of the set of
numbers in that row is $1$. One computes $f(3)=2$; $f(4)=4$; and so on.
Computation of the exact value of $f(k)$ is slightly tricky and is not
reported here. For all $k$, $f(k)\geq2$. To summarize, for each fixed
row there are a.s. ultimately $f(k)k!$ disjoint (open) invariant sets
following row normalization, making for $[f(k)k!]^n$ a.s. ultimately
(open) invariant sets simultaneously for all n rows. If we count a.s.
ultimately invariant sets subsequent to column normalization, then
entirely analogous arguments result in a.s. ultimately $[f(n)n!]^k$
disjoint (open) sets simultaneously for all columns.

From computations, after a row normalization the surface area of the
sphere in $k$-space orthogonal to the equiangular line---that
corresponds to only one row of $\mathbf{X}^{(m)}$---has surface area
$\approx\sqrt{ (\frac{2}{e} )}  (\frac{2\pi e}{k-3}
) < 1$ for $k \geq21$. The expression $\rightarrow0$ as $k\nearrow
\infty$. Even for $k=4$, the quantity is only about 14.7 (larger than
the actual value). Remember that there are at most $[f(k)k!]$ a.s.
ultimately nonempty ``invariant sets'' for row normalization. Thus one
sees that for $k$ large the quantity $ (\frac{\mathrm{surface\ area\ }\operatorname{Sph}(k)}{|\mathrm{invariant\ sets\ of\ }\operatorname{Sph}(k)|} )^{n}$
is nearly $0$.

\section{Computational results and applications}

We include three examples to highlight and illustrate some
computational aspects of our iterative procedure. The first two
examples are studies by simulation whereas the third example is an
implementation on a real dataset.

For the simulation study, we consider a 3-by-3 matrix and a 10-by-10
matrix both with entries generated independently from a uniform
distribution on [0, 1]. For a given matrix, the algorithm
computes/calculates the following 4 steps at each iteration:
\begin{longlist}[(a)]
\item[(a)] mean polish the column,
\item[(b)] stand deviation polish the column,
\item[(c)] mean polish the row and
\item[(d)] stand deviation polish the row.
\end{longlist}

These fours steps, which constitute one iteration, are repeated
until\break
``converg\-ence''---which we define as when the difference in some
norm-squared (the quadratic/Frobenius norm in our case) between two
consecutive iterations is less than some small prescribed value---which
we take to be 0.00000001 or $10^{-8}$.

%
%
%
%
%
%
%
%

\subsection{Example 2: Simulation study on a 10-by-10 dimensional example}

We proceed now to illustrate the convergence of the successive row and
column mean--standard deviation polishing for the simple 10-by-10
dimensional example cited. The algorithm took 15 iterations to
converge. The initial matrix, the final solution, and relative (and log
relative) difference for the 15 iterations follow:
%
\begin{eqnarray}
\hspace*{50pt}X^{0} &=&
\left(\matrix{
0.8145 & 0.3551 & 0.7258 & 0.3736 & 0.0216    \cr
0.7891 & 0.9970 & 0.3704 & 0.0875 & 0.9106    \cr
0.8523 & 0.2242 & 0.8416 & 0.6401 & 0.8006    \cr
0.5056 & 0.6525 & 0.7342 & 0.1806 & 0.7458    \cr
0.6357 & 0.6050 & 0.5710 & 0.0451 & 0.8131    \cr
0.9509 & 0.3872 & 0.1769 & 0.7232 & 0.3833    \cr
0.4440 & 0.1422 & 0.9574 & 0.3474 & 0.6173    \cr
0.0600 & 0.0251 & 0.2653 & 0.6606 & 0.5755    \cr
0.8667 & 0.4211 & 0.9246 & 0.3839 & 0.5301    \cr
0.6312 & 0.1841 & 0.2238 & 0.6273 & 0.2751}\right.\nonumber
\\[-8pt]\\[-8pt]
&&\hspace*{5pt}{}\left.\matrix{0.2486 & 0.0669 & 0.2178 & 0.6766 & 0.6026\cr
 0.4516 & 0.9394 & 0.1821 & 0.9883 & 0.7505\cr
 0.2277 & 0.0182 & 0.0418 & 0.7668 & 0.5835\cr
 0.8044 & 0.6838 & 0.1069 & 0.3367 & 0.5518\cr
 0.9861 & 0.7837 & 0.6164 & 0.6624 & 0.5836\cr
 0.0300 & 0.5341 & 0.9397 & 0.2442 & 0.5118\cr
 0.5357 & 0.8854 & 0.3545 & 0.2955 & 0.0826\cr
 0.0871 & 0.8990 & 0.4106 & 0.6802 & 0.7196\cr
 0.8021 & 0.6259 & 0.9843 & 0.5278 & 0.9962\cr
 0.9891 & 0.1379 & 0.9456 & 0.4116 & 0.3545
}\right),\nonumber
%
%
\\
%
X^{\mathit{final}} &=&
\left(\matrix{
1.2075 & 0.2139 & 0.8939 & 0.2661 & -2.0026       \cr
-0.0736 & 1.7222 & -1.2202 & -1.0461 & 0.6465     \cr
0.8858 & -0.8659 & 0.8816 & 0.7930 & 0.9515       \cr
-0.9296 & 1.5223 & 0.6537 & -0.7661 & 0.9476      \cr
-0.8358 & 0.8041 & -0.7288 & -2.0057 & 1.0328     \cr
1.4926 & 0.1374 & -1.2120 & 1.1351 & -0.5035      \cr
-0.5156 & -0.7494 & 1.5647 & 0.2025 & 0.5610      \cr
-1.8680 & -1.1055 & -0.6428 & 0.9269 & 0.3515     \cr
0.3596 & -1.0158 & 0.8070 & -0.5547 & -1.1339     \cr
0.2771 & -0.6632 & -0.9973 & 1.0490 & -0.8509}\right.\nonumber
\\[-8pt]\\[-8pt]
&&\hspace*{5pt}{}\left.\matrix{ -0.5881 & -1.2477 & -0.4157 & 1.1023 & 0.5705   \cr
 -0.8172 & 0.5144 & -1.1740 & 1.3022 & 0.1458                  \cr
-0.9498 & -1.6621 & -1.1469 & 1.0831 & 0.0298                  \cr
0.9361 & 0.5467 & -1.3402 & -1.3775 & -0.1931                  \cr
 1.4824 & 0.5929 & 0.0202 & 0.2768 & -0.6390                   \cr
-1.2741 & 0.1766 & 1.3125 & -1.1642 & -0.1005                  \cr
 0.2646 & 1.3840 & -0.0059 & -0.7521 & -1.9537                 \cr
 -0.8323 & 1.2448 & 0.1167 & 0.8827 & 0.9259                   \cr
 0.1669 & -0.5895 & 1.0581 & -1.0805 & 1.9828                  \cr
 1.6114 & -0.9601 & 1.5752 & -0.2727 & -0.7685}%
\right),\nonumber
%
\\
%
&&\hspace*{-33pt}\mbox{Successive Difference}\nonumber
\\[-8pt]\\[-8pt]
&&\hspace*{-33pt}\qquad=\pmatrix{
\mbox{Iteration no.} & \mbox{Difference} & \mbox{log(difference)}\cr
1 & 84.1592 & 4.4327 \cr
2 & 1.2860 & 0.2516 \cr
3 & 0.1013 & -2.2897 \cr
4 & 0.0144 & -4.2402 \cr
5 & 0.0029 & -5.8434 \cr
6 & 0.0007 & -7.2915 \cr
7 & 0.0002 & -8.6805 \cr
8 & 0.0000 & -10.0456 \cr
9 & 0.0000 & -11.4000 \cr
10 & 0.0000 & -12.7492 \cr
11 & 0.0000 & -14.0955 \cr
12 & 0.0000 & -15.4403 \cr
13 & 0.0000 & -16.7841 \cr
14 & 0.0000 & -18.1272 \cr
15 & 0.0000 & -19.4699
}\nonumber
.
\end{eqnarray}

\begin{figure}[b]

\includegraphics{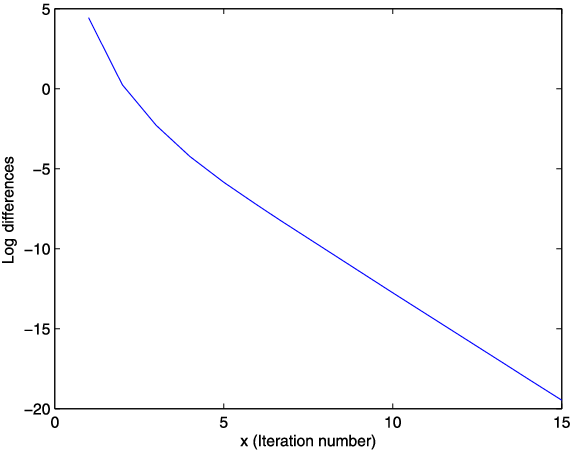}

\caption{Relative differences at each iteration on the log
scale---10-by-10 dimensional example.}\label{fig:51}
\end{figure}

\noindent We note once more how the relative differences decrease
linearly on the log scale (though empirically) and how this is again
suggestive of the rate of convergence. As both the figure (see Figure~\ref{fig:51}) and the vector of relative differences indicate, there is a
substantial jump at iteration 2, and then the curve behaves linearly.

\begin{figure}[b]

\includegraphics{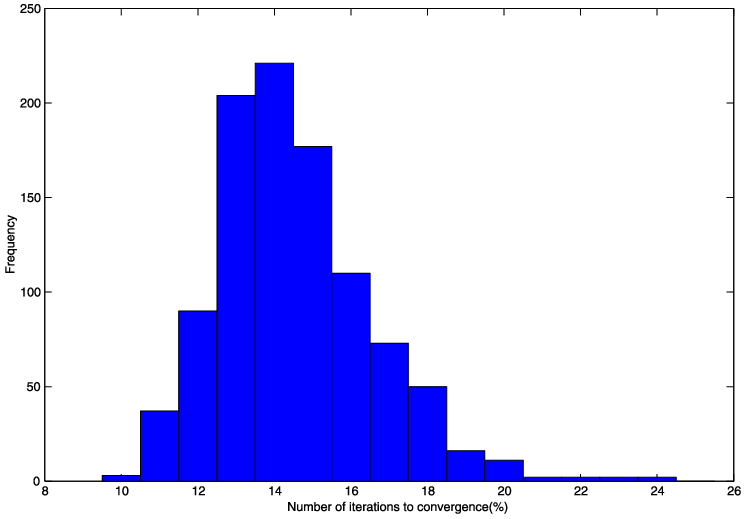}

\caption{Distribution of the number of steps to convergence.}\label{fig:52}
\end{figure}

The whole procedure takes about 0.37 seconds on a standard
modern laptop computer and terminates after 15 iterations. It might
appear that the increase in the number of iterations increases with
increase in dimension. For instance, the number of iterations goes from
9 to 15 as we go from dimension 3 to 10. We should, however, bear in
mind that when we go from dimension 3 to 10 the ``tolerance level'' is
kept constant at 0.00000001. The number of elements that must be close
to their respective limiting values, however, goes from 9 in the
3-dimensional case and to 100 in the 10-dimensional case. The rapidity
of convergence was explored further, and the process above was repeated
over $1000$ simulations. The convergence proves to be stable in the
sense that the mean and standard deviation of the number of steps until
convergence over the $1000$ simulations are 14.5230 and 2.0331,
respectively. A histogram of the number of steps till convergence is
given below (Figure~\ref{fig:52}).

%

A closer look at the vector of successive differences
suggests that the ``bulk of the convergence'' is achieved during the
first iteration. This seems reasonable since the first steps render the
resulting columns, then the rows and the members of their respective
unit spheres (and even within cited subsets of them). Convergence then
is only within these spheres. Our numerical results also indicate the
the sequence of matrices $X^{(i)}$ changes most drastically during this
first iteration. It suggests that mean polishing to a larger extent is
responsible for the rapidity of convergence and is reminiscent of the
result in Lemma~\ref{lm:lemma1.3.1} which states that if only row and column mean
polishing are performed, then convergence is achieved immediately. We
explore this issue further by looking at the distance between $X^{(1)}$
and $X^{(\mathit{final})}$ and compare it to the distance between $X^{(0)}$ and
$X^{(\mathit{final})}$. The ratio of these two distances follows:
\[
\mathit{Ratio} = \frac{\operatorname{dist}(X^{(1)}, X^{(\mathit{final})})}{\operatorname{dist}(X^{(0)},
X^{(\mathit{final})})}.
\]

\begin{figure}[b]

\includegraphics{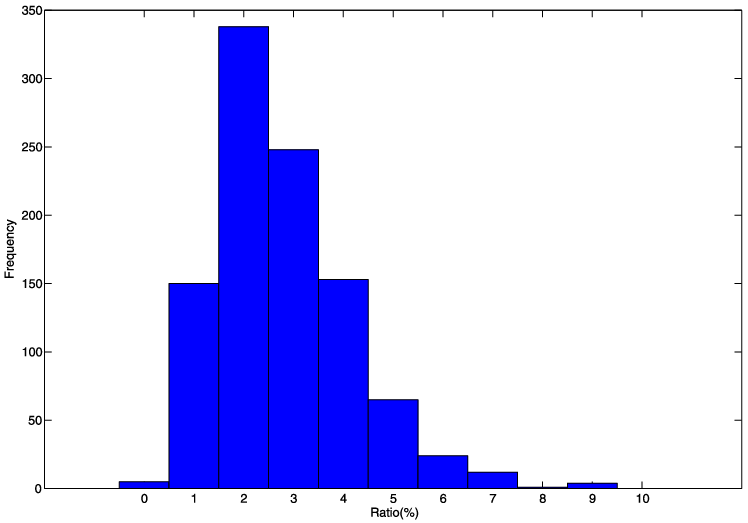}

\caption{Distribution of distance to limit after 1-step as a proportion
of distance to limit from initial-step.}\label{fig:53}
\end{figure}

\noindent For our 10-by-10 example, we simulated 1000 initial starting
values and implemented our successive normalization procedure. The
average value of the distance from the first iterate to the limit as a
proportion of the total distance to the limit from the starting value
is only 2.78$\%$. One could interpret this heuristically as saying that
on average the crucial first step does as much as 97.2$\%$ of the work
towards convergence. We therefore confirm that the bulk of the
convergence is indeed achieved in the first step (termed as a
``one-step analysis'' from now onwards). The distribution of the ratio
defined above is graphed in the histogram below (Figure~\ref{fig:53}). We
also note that none of the 1000 simulations yielded a ratio of over
10$\%$.\looseness=1

%

Yet a another illuminating perspective of our successive
normalization technique is obtained when we track the number of sign
changes in the individual entries of the matrices from one iteration to
the next. Please remember that this is related to the ``invariant
sets'' that were described in Section~\ref{domains_of_attraction}.
Naturally, one would expect the vast majority of sign changes to occur
in the first step as the bulk of the convergence is achieved during
this first step. We record the number of sign changes at each
iteration, as a proportion of the total number of sign changes until
convergence, over 1000 simulations, in our 10-by-10 case. The results
are illustrated in the table below\footnote{Since the number of
iterations to convergence depends on the starting point, the length of
the vector of the number of sign changes will vary accordingly. We
summarize this vector by averaging over all the 1000 simulations the
relative frequency of the number of sign changes for the first nine
iterations. The first nine iterations were chosen as each of the 1000
simulations required at least 9 iterations to converge.} (see Table~\ref{table:51}). An empirical study of the occurrence of the sign changes
reveals interesting heuristics. We note that on average 95\% of sign
changes occur during the first step and an additional 3\% in the next
step. The table also demonstrates that as much as 99\% of sign changes
occur during the first three iterations. When we examine infrequent
cases where there is a sign change well after the first few iterations,
we observe that the corresponding limiting value is close to zero thus
indicating that a sign change well into the successive normalization
technique (i.e., a change from positive to negative or vice versa)
amounts to a very small change in the actual magnitude of the
corresponding value of the matrix.

%
\begin{table}
\tablewidth=10cm
\caption{Distribution of the occurrence of sign changes}\label{table:51}
\begin{tabular*}{\tablewidth}{@{\extracolsep{\fill}}lcc@{}}
\hline
\textbf{Iteration no.} & \textbf{Relative frequency} & \textbf{Relative cumulative frequency} \\
\hline
\phantom{1}1 & \phantom{1}94.97\% & \phantom{1}94.97\% \\
\phantom{1}2 & \phantom{19}3.15\% & \phantom{1}98.12\% \\
\phantom{1}3 & \phantom{19}1.03\% & \phantom{1}99.15\% \\
\phantom{1}4 & \phantom{19}0.40\% & \phantom{1}99.55\% \\
\phantom{1}5 & \phantom{19}0.20\% & \phantom{1}99.75\% \\
\phantom{1}6 & \phantom{19}0.12\% & \phantom{1}99.87\% \\
\phantom{1}7 & \phantom{19}0.05\% & \phantom{1}99.92\% \\
\phantom{1}8 & \phantom{19}0.03\% & \phantom{1}99.95\% \\
\phantom{1}9 & \phantom{19}0.02\% & \phantom{1}99.97\% \\
10\mbox{ and above} & \phantom{19}0.03\% & 100.00\% \\[3pt]
Total & 100.00\% & -- \\
\hline
\end{tabular*}
\end{table}

We conclude this example by investigating more thoroughly whether the
dimensions of the matrices have an impact on either the rapidity of
convergence and/or on the one-step analysis. The following table gives
the mean and standard deviations of the number of iterations needed for
convergence for various values of the dimension of the matrix, denoted
by $p$ and $n$ when keeping the total number of cells in the matrix
constant.\footnote{Or approximately constant.} Once more our successive
normalization procedure is applied to $1000$ uniform random starting
values. Results of this exercise are given in the table below (see
Table~\ref{table:52}).

%
\begin{table}
\tablewidth=8.5cm
\caption{Rapidity of convergence and one-step analysis for various\break $p$
and $n$ combinations}\label{table:52}
\begin{tabular*}{\tablewidth}{@{\extracolsep{\fill}}l c c c c@{}}
\hline
$\bolds{p}$ & \textbf{3} & \textbf{4} & \textbf{5} & \textbf{10} \\
\hline
$\bolds{n}$ & \textbf{33} & \textbf{25} & \textbf{20} & \textbf{10} \\
\hline
$\operatorname{mean}$(count)            &34.0870             & 22.3500           & 18.0720           & 14.5790 \\
$\operatorname{std}$(count)             &\phantom{0}5.3870   & \phantom{0}3.3162 & \phantom{0}2.5477 & \phantom{0}2.1099
\\[3pt]
$\operatorname{mean}$(ratio)            &\phantom{0}2.7026   & \phantom{0}2.5812 & \phantom{0}2.6237 & \phantom{0}2.7207 \\
$\operatorname{std}$(ratio)             &\phantom{0}1.8483   & \phantom{0}1.5455 & \phantom{0}1.4768 & \phantom{0}1.2699 \\
\hline
\end{tabular*}
\end{table}

We find that when $n$ and $p$ are close, convergence appears
to be faster, keeping everything else constant. Interestingly enough, a
one-step analysis performed for the different scenarios above tends to
suggest that the one-step ratio, defined as the distance from the first
iterate to the limit as a proportion of the total distance to the limit
from the starting value, seems largely unaffected by the row or column
dimension of the problem.

\subsection{Example 3: Simulation study on a 5-by-5 dimensional example}

We now proceed to further investigate the successive normalization
procedure when one begins with column mean-standard deviation polishing
followed by row mean-standard deviation polishing or vice versa on a
simple 5-by-5 dimensional example. The theory developed in the previous
sections proves convergence of the successive normalization procedure
whether the first normalization that is performed on the matrix is row
polishing or column polishing.

\begin{figure}

\includegraphics{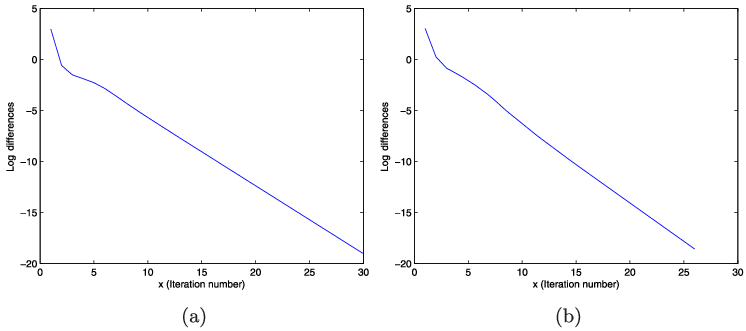}

\caption{Relative differences at each iteration on the log scale for
5-by-5 dimensional example~\textup{(a)}~starting with column
polishing \textup{(b)} starting with row polishing.}\label{fig:54}
\end{figure}

The algorithm took 30 iterations to converge when one begins with
column mean-standard deviation polishing, and when one begins with row
mean-standard deviation polishing it took 26 iterations to converge.
The initial matrix, the final solutions, log relative differences and
their respective plots for both approaches are given below (see
Figure~\ref{fig:54}):

%
\begin{eqnarray}
&&\hspace*{6pt}X^{0} =
\pmatrix{
0.6565 & 0.2866 & 0.7095 & 0.4409 & 0.8645 \cr
0.3099 & 0.3548 & 0.9052 & 0.8758 & 0.0210 \cr
0.3316 & 0.5358 & 0.8658 & 0.8650 & 0.0768 \cr
0.1882 & 0.9908 & 0.1192 & 0.3552 & 0.3767 \cr
0.1007 & 0.0282 & 0.9553 & 0.6311 & 0.1492
}
,
%
\\
%
&&X^{\mathit{final}}_{\mathrm{starting\ with\ column\ polishing}}\nonumber
\\[-8pt]\\[-8pt]
&&\qquad=
\pmatrix{
1.6360 & -0.4320 & -0.7863 & -1.0548 & 0.6371 \cr
0.1093 & -1.2446 & 0.9477 & 1.2170 & -1.0295 \cr
-0.6979 & 1.1193 & -0.1716 & 1.1399 & -1.3897 \cr
0.2748 & 1.2421 & -1.3091 & -1.0112 & 0.8034 \cr
-1.3223 & -0.6848 & 1.3192 & -0.2907 & 0.9786
}\nonumber
,
%
\\
%
&&X^{\mathit{final}}_{\mathrm{starting\ with\ row\ polishing}}\nonumber
\\[-8pt]\\[-8pt]
&&\qquad=
\pmatrix{
1.4956 & -0.4243 & -0.7386 & -1.1620 & 0.8293 \cr
0.3816 & -0.9267 & 0.5915 & 1.3267 & -1.3731 \cr
-1.2158 & 1.1775 & 0.1052 & 0.9966 & -1.0634 \cr
0.3478 & 1.2181 & -1.4096 & -0.9138 & 0.7573 \cr
-1.0092 & -1.0446 & 1.4514 & -0.2475 & 0.8499
}\nonumber
,
%
\\
%
&&\mbox{Successive Difference}\nonumber
\\[-8pt]\\[-8pt]
&&\qquad=
\pmatrix{
 & \mbox{Difference} & \mbox{log(difference)}\cr
 & \mbox{starting with} & \mbox{starting with} \cr
\mbox{Iteration no.} & \mbox{column polishing} & \mbox{row polishing} \cr
1 &2.9646 & 3.0255 \cr
2 & -0.5858 & 0.2539 \cr
3 & -1.5082 & -0.8731 \cr
4 &-1.8814 & -1.4650 \cr
\ldots & \ldots & \ldots \cr
\ldots & \ldots & \ldots \vspace*{4pt}\cr
24 & -15.0375 & -17.0730 \cr
25 & -15.7028 & -17.8229 \cr
26 & -16.3679 & -18.5728 \cr
27 & -17.0331 & - \cr
28 & -17.6983 & - \cr
29 & -18.3635 & - \cr
30 & -19.0287 & -
}\nonumber
.
\end{eqnarray}

%

As expected, the final solutions are different. The simulations were
repeated with different initial values, and we note that the
convergence patterns (as illustrated in Figure~\ref{fig:54}) are similar
whether the procedure starts with column polishing or row polishing,
though the actual number of iterations required to converge can vary.

\subsection{Example 4: Gene expression data set: 20426-by-63
dimensional example}

We now illustrate the convergence of the successive row and column
mean--standard deviation polishing for a real-life gene expression
example, a 20426-by-63 dimensional example. This dataset arose
originally from a study of human in-stent restenosis by Ashley et al.
\cite{ashley}. The algorithm took considerably longer in terms of time
and computer resources but converged in eight iterations. The initial
matrix and the final are too large to display, but the relative (and log
relative) differences for the eight iterations are given subsequently:

%
\begin{eqnarray}
&&\mbox{Successive Difference}\nonumber
\\[-8pt]\\[-8pt]
&&\qquad= 1.0e\!+\!005 *
\pmatrix{
\mbox{Iteration no.} & \mbox{Difference} & \log\mbox{(difference)}\cr
1 & 1.0465 & 11.5583 \cr
2 & 0.0008 & 4.4030 \cr
3 & 0.0000 & -0.2333 \cr
4 & 0.0000 & -4.7582 \cr
5 & 0.0000 & -9.2495 \cr
6 & 0.0000 & -13.7130 \cr
7 & 0.0000 & -18.1526 \cr
8 & 0.0000 & -22.5717
}\nonumber
.
\end{eqnarray}

%
%
%
%
%
%

%
\begin{figure}

\includegraphics{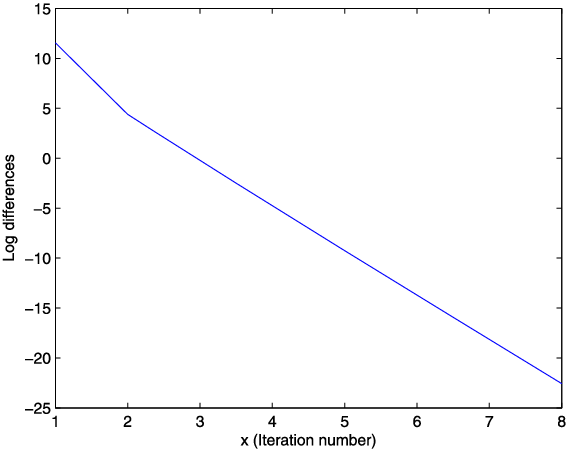}

\caption{Relative differences at each iteration on the log
scale---20426-by-63 dimensional example.}\label{fig:55}
\end{figure}

\noindent Note once more how the relative differences decrease
linearly on the log scale (though empirically) and is once again
suggestive of the rate of convergence. As both the figure (see Figure~\ref{fig:55}) and the vector of relative differences indicates, there is a
jump between iteration 1 and 2 and then the curve behaves linearly.

 Additionally the whole procedure takes about 853.2 seconds or
approximately 14.22 minutes on a desktop computer\footnote{2~GHz Core 2
Duo processor and 2~GB of RAM.} versus 0.4 seconds for the 10-by-10
example. However, the algorithm terminates after only eight iterations.
In this example the number of iterations does NOT change with the
increase in dimensionality. It may make sense to investigate this
behavior more thoroughly, empirically, using simulation for rectangular
but not square matrices. It seems that the ratio of the two dimensions
or the minimum of the two dimensions may play a role. We should also
bear in mind that the tolerance level, which is the sum of the
individual differences squared, has been kept constant at 0.00000001.

\section{Conclusion}

In this section we attempt to lend perspective to our results and to
point the way for future developments. Readers please note that for
rectangular $n\times k$ arrays of real numbers with $\min(k,n)\geq3$,
the technique beginning with rows (alternatively columns) and
successively subtracting row (column) means and dividing resulting
differences, respectively, by row (column) standard deviations
converges for a subset of Euclidean $\mathcal{R}^{n \times k}$ whose
complement has Lebesgue measure 0. The limit is row and column
exchangeable given the Gaussian probability mechanism that applies in
our theoretical arguments. We do not offer other information on the
nature of the exact set on which successive iterates converge. A single
``iteration'' of the process we study has four steps, two each,
respectively, for rows and columns. Note that on the set for which the
algorithm converges, convergence seems remarkably rapid, exponential or
even faster, perhaps because after a half an iteration, the rows
(alternately columns) lie as $n$ (respectively $k$) points on the
surface of the product of relevant unit spheres. Further iterations
adjust only directions, not lengths.

Viewing the squares of the entries as the terms of a backwards
martingale shows maximal inequalities for them, and therefore
implicitly contains information on ``rates of convergence'' of the
squares; but these easy results appear far from the best one might
establish. Our arguments for (almost everywhere) convergence of the
original signed entries do not have information regarding rates of
convergence. One argues easily that if successive iterates converge,
and no limiting entry is $0$, then after finitely many steps (the
number depending on the original values and the limiting values), signs
are unchanged. In our examples of small dimension, evidence of this can
be made explicit. In particular we observe empirically that the vast
majority of sign changes that are observed do indeed take place in the
first few iterations. Any sign changes that are observed well after the
first few iterations correspond to sign changes around entries with
limiting values close to zero. We also have no information on
optimality in any sense of the iterated transformations we study. One
reason for our thinking that our topic is inherently difficult is that
we were unable to view successive iterates as ``contractions'' in any
sense familiar to us.

If we take any original set of numbers, and multiply each number by the
realized value of a positive random variable with arbitrarily heavy
tails, then convergence is unchanged. Normalization requires that after
half a single iteration the same points on the surface of the relevant
unit spheres are attained, no matter the multiple. The message is that
what matters for convergence are the distributions induced on the
surfaces of spheres after each half iteration, and not the otherwise
common heaviness of the tails of the probability distributions of
individual entries.

\section*{Acknowledgments}
The authors gratefully acknowledge Bradley Efron for introducing them
to the research question addressed in the paper.
They also thank Johann Won and Thomas Quertermous for useful
discussions and Thomas Quertermous for granting
us access to and understanding of data we study and reference; Bonnie
Chung and Cindy Kirby are also
acknowledged for administrative assistance.

\printaddresses

\end{document}


\begin{frontmatter}

\title{Correction\vspace*{6.6pt}\\
Successive normalization of rectangular arrays}
\runtitle{Correction}

\pdftitle{Correction:
Successive normalization of rectangular arrays}

\textit{Ann. Statist.} \textbf{38} (2010) 1638--1664

\begin{aug}
\author{\fnms{Richard A.} \snm{Olshen}\corref{}\ead[label=e1]{olshen@stanford.edu}}
\and
\author{\fnms{Bala}~\snm{Brajaratnam}\ead[label=e3]{brajarat@stanford.edu}}
\runauthor{R. A. Olshen and B. Brajaratnam}
\affiliation{Stanford University}
\address{Departments of Health Research and Policy\\
\quad Electrical Engineering, and Statistics\\
Stanford University\\
Stanford, California 94305-4065\\
USA\\
\printead{e1}}

\address{Departments of Statistics\\
\quad and of Environmental\\
\quad Earth System Science\\
Stanford University\\
Stanford, California 94305-4065\\
USA\\
\printead{e3}}
\end{aug}

\received{\smonth{4} \syear{2013}}
\revised{\smonth{9} \syear{2013}}

%
\begin{abstract}
We report an error in Theorem~4.1 of the paper identified in the title.
That theorem is not true in general. However, the main mathematical
message in the paper as well as all subsidiary mathematical remarks
and all simulations and numerical results are correct. As the title
suggests, the paper concerns successive normalization of rectangular
arrays. We learned the algorithm from Bradley Efron. Here we sketch
the algorithm, state the mistake in our paper, sketch our approach
to proving its main theorem, and give reference to a paper by the
authors where many more details can be found.
\end{abstract}

%
\begin{keyword}[class=AMS]
\kwd{62H05}
\kwd{60F15}
\kwd{60G46}
\end{keyword}
\begin{keyword}
\kwd{Normalization}
\kwd{standardization}
\end{keyword}

\end{frontmatter}

For a given matrix, Efron's algorithm involves the following four
successive steps at each iteration:

\begin{enumerate}
\item Mean polish each of the $J$ columns.
\item Standard deviation polish each column.
\item Mean polish each of the $I$ rows.
\item Standard deviation polish each row.
\end{enumerate}

These fours steps, which constitute one iteration, are repeated until
``convergence.''

Theorem~4.1 of the paper by Olshen and Rajaratnam \cite{2} is false as it
stands. Writing ${\cal G}^{(i)}_{2m-1}$ correctly does not change
things. One reason is that, in the notation of the paper, typically
\[
E\bigl\{\bigl( X^{(2m-1)}_{i \pi(1)} \bigr)^2 I_B
I_Q \bigr\}\quad \mbox{is not}\quad
E\bigl\{\bigl( X^{(1)}_{i \pi(1)} \bigr)^2 I_B
I_Q \bigr\}.
\]

Efforts toward finding a nontrivial sub-$\sigma$-field of ${\cal
G}^{(i)}_{2m-1}$ for which there might be equality failed. However, the
goal of the theorem does obtain; but it is a corollary,\vadjust{\goodbreak} not a tool, in
proving the a.s. convergence of Efron's algorithm applied to
$\mathbf{X}_{n \times k} \sim N_{n \times k} (\mathbf{0}, \mathbf{I})$.
That is, if, as is the case, the argument from the bottom of page 1648
through the end of Section~4 (of \cite{2}) can be proven without leaning on
Theorem~4.1, then not only does the purpose for Theorem~4.1 hold, but
in fact much more. Namely, write
${\cal F}_i^{(m)}$ for the $\sigma$-field of the first display on page
1647. For each fixed $i$, these $\sigma$-fields are nested (decreasing)
in $m$. Furthermore, off the set of probability $0$ implicit in (9),
page 1646, $\sum_{i,j} (X_{ij}^{(q)})^2 = IJ$ for all $(i, j, q)$. We
assume henceforth that $\mathbf{x}$ lies outside this ``cursed'' set.
Therefore, for any real-valued, continuous function $f$ on
${\cal R}^{n\times k}, \lim_{m\rightarrow\infty} E\{ f (X_{ij}^{(q)} )
| {\cal F}_i^{(m)}\}$ exists almost surely for every fixed $q$ (as $m$
increases without bound), but also, as a consequence of Theorem~2 of
the paper by Blackwell and Dubins \cite{1}, the same holds with $q$ replaced
by $m$.\looseness=-1

Now we sketch an argument by which one shows directly that
simultaneously in $(i,j), X_{ij}^{(m)}$ converges a.s. as $m$ increases
without bound. First, we study $(S_j^{(2m-1)})^2$ as in the last
display on page 1648, and we sketch the argument that
$P\{ \overline{\lim}_m (S_j^{(2m-1)})^2 > 0 \} = 1$. As in \cite{2}, let
$A= \{ \overline{\lim}_m (S_j^{(2m-1)})^2 = 0\}$.
For a more detailed argument, please see \cite{3}.

Since the entries of $\mathbf{X}$ are independent, $\mathbf{X}$ is row
and column exchangeable. This property is inherited by $\mathbf
{X}^{(q)}$ for every $q$. Because all entries (for $q \ge1$) are
bounded, $E\{ X_{ij}^{(q)}\}$ and $E\{ (X_{ij}^{(q)})^2 \}$ exist and
are finite (with fixed bound that applies to all $q$). Exchangeability
implies that all $E\{ X_{ij}^{(q)} \} = 0$, and all $E\{
(X_{ij}^{(q)})^2 \} = 1$. Bounded convergence implies that if
$(S_j^{(2m-1)})^2$ tends to $0$ along a subsequence as $m$ increases,
not only is the limit bounded as a function of $\mathbf{x}$, but also
the limit random variable has expectation $0$. Necessarily every almost
sure subsequential limit in $m$ of the random variables $\overline
{X}_{\cdot j}^{(2m-1)}$ has mean $0$. Likewise, every almost sure
subsequential limit in $m$ of the random variables
$(X_{ij}^{(2m-1)})^2$ has expectation 1. All are bounded as functions
of $\mathbf{x}$. One consequence of these things is that $P(A)=0$.
Please note that the statements about $\underline{\lim}_m
X_{ij}^{(2m-1)}$ and $\overline{\lim}_m X_{ij}^{(2m-1)}$ are wrong.

Continue to the paragraph on page 1649 that begins, ``Again, let\ldots''
Define $m_q$ as in the paper. The cardinality of $E$ is at least 2. (In
fact, the cardinality of $E$ exceeds 2.) Now proceed to the next
paragraph. We define the last three displays on page 1649 slightly
differently; thus, we require that along $m_q$,
\begin{eqnarray*}
&&\lim_{m_q} X_{i_o j}^{(2m_q -1)}\quad \mbox{and}\quad \lim
_{m_q} X_{i_0
j}^{(2m_q)} \qquad\mbox{both exist;}
\\
&&\lim_{m_q} X_{i_1 j}^{(2m_q -1)}\quad \mbox{and}\quad \lim
_{m_q} X_{i_1
j}^{(2m_q)}\qquad \mbox{both exist};\quad \mbox{and}
\\
&&\lim_{m_q} \bigl| X_{i_o j}^{(2m_q)} -
X_{i_1 j}^{(2m_q)} \bigr| \qquad\mbox{exists and is positive}.
\end{eqnarray*}

If there were no such subsequence, our proof would be complete (though
previous arguments show that this is not a possibility). Proceed to the
first two displays on page 1650. The argument for $S_j^{(2m-1)}$
tending to 1 holds not because the cited difference tends to $0$
(something we do not know yet), but instead because of the third
display. The remainder stays the same. Theorem~\ref{th4.2} follows.

\renewcommand{\thetheorem}{4.2}
\begin{theorem}\label{th4.2}
So long as I and J are at least $3$, Efron's algorithm converges almost
surely for $\mathbf{X}$, and therefore on a Lebesgue set of entries
with complement a set of Lebesgue measure $0$.
\end{theorem}

%

\printaddresses